\documentclass[a4paper, 10pt, twoside]{article}

\usepackage{amsmath, amscd, amsfonts, amssymb, amsthm, latexsym, url, color, todonotes, bm, framed, rotating, enumerate} 
\usepackage{graphicx}
\usepackage[left=1in, right=1in, top=1.2in, bottom=1in, includefoot, headheight=13.6pt]{geometry}
\usepackage{mathtools}
\usepackage[all]{xy}
\usepackage{caption} %

\usepackage{marvosym}
\usepackage{tikz}
\usetikzlibrary{babel}
\usepackage{tikz-cd}
\usepackage{thmtools}
\usepackage{thm-restate}
\newcommand{\authormark}[1]{\textsuperscript{#1}}

\newcommand{\eliaMark}{*}

\newlength{\sprRuleWidth}
\newcommand{\sprrule}{\rule{\sprRuleWidth}{0.3pt}}

\newcommand{\corrsym}{\raisebox{0.15ex}{\small\Letter}}

\makeatletter
\renewcommand{\footnoterule}{%
  \kern-3pt\sprrule\kern 2.6pt%
}
\makeatother

\newlength{\corrindent}
\usepackage{tablefootnote}
\usepackage[section]{placeins}

\usepackage{threeparttable}
\usepackage{diagbox}

\usepackage[colorlinks,breaklinks=true]{hyperref}
\usepackage[figure,table]{hypcap}
\hypersetup{
	bookmarksnumbered,
	pdfstartview={FitH},
	citecolor={black},
	linkcolor={black},
	urlcolor={black},
	pdfpagemode={UseOutlines}
}
\makeatletter
\newcommand\org@hypertarget{}
\let\org@hypertarget\hypertarget
\renewcommand\hypertarget[2]{%
  \Hy@raisedlink{\org@hypertarget{#1}{}}#2%
} 
\makeatother 

\newtheorem{theorem}{Theorem}[section]
\newtheorem{lemma}[theorem]{Lemma}

\newtheorem{proposition}[theorem]{Proposition}

\theoremstyle{definition}
\newtheorem{definition}[theorem]{Definition}

\newtheorem{remark}[theorem]{Remark}

\newcommand{\xysquare}[8]{
\[\xymatrix{
#1 \ar@{#5}[r] \ar@{#6}[d] & #2 \ar@{#7}[d]\\
#3 \ar@{#8}[r] & #4
}\]
}

\newcommand{\bb}{\mathbb}

\newcommand{\comment}[1]{}

\renewcommand{\phi}{\varphi}

\usepackage{calc} 
\usepackage{enumitem}

\newlist{caselist}{enumerate}{1}

\setlist[caselist]{%
  label=\textnormal{Case (\roman*)}, 
  align=left,                        
  leftmargin=*,                      
  itemsep=0pt                        
}

\renewcommand{\cal}{\mathcal}
\renewcommand{\hat}{\widehat}

\DeclareMathOperator{\Spec}{Spec}

\newcommand{\rcoeq}[3]{\xymatrix{ #1\ar@/^3mm/[r]^f \ar@/_3mm/[r]_g & #2 \ar[l]_e\ar[r] & #3}}

\usepackage[bbgreekl]{mathbbol}
\DeclareSymbolFontAlphabet{\mathbbm}{bbold}

\usepackage{fancyhdr}

\usepackage{sectsty}
\sectionfont{\Large\sc\centering}
\partfont{\centering}

\begin{document}

\title{Hypersurfaces that are far from ruled}

\author{%
  Elia Fiammengo\authormark{\eliaMark}%
}

\date{}

\maketitle

\begingroup
\makeatletter
\renewcommand{\@makefntext}[1]{\noindent#1}
\makeatother

\renewcommand\thefootnote{}
\footnotetext{%
  \footnotesize

  \noindent\sprrule\par
  \vspace{0.35\baselineskip}

  \noindent\corrsym\hspace{0.4em}Elia Fiammengo\par
  \noindent\hspace*{\corrindent}\texttt{elia.fiammengo@stud.uni-heidelberg.de}\par
  \vspace{0.35\baselineskip}

  \noindent\authormark{\eliaMark}\; Universit\"at Heidelberg, Mathematisches Institut,
  Im Neuenheimer Feld~205, 69120 Heidelberg, Germany\par
}%
\endgroup

\setcounter{footnote}{0}

\begin{abstract}

We show that a very general complex hypersurface of degree $d$ in
$\bb P^{n+1}$ is not ruled whenever $n\geq3$ and
$d\geq2\lceil(n+2)/3\rceil$. This improves Kollár's original degree bound to the bound subsequently obtained by Totaro for stable irrationality and, in particular, includes very general quartic fourfolds. 
More precisely, we find new bounds for the minimal degree of a dominant rational map from a very general hypersurface to a ruled variety or a conic bundle, improving on previous bounds of Chen--Stapleton.
The proofs combine Kollár's reduction to positive characteristic
with the Jacobson correspondence for inseparable extensions.

\end{abstract}
\begingroup

\smallskip

\par
\endgroup

\section{Introduction}
A variety $X$ of dimension $n$ is ruled if there exists a birational map $X\dashrightarrow  Y\times \bb P^1$ for some variety $Y$ of dimension $n-1$. 
It is an interesting question to determine if a very general degree $d$ complex hypersurface $X_d\subset \bb P^{n+1}$ is ruled and it is expected that this is not the case for $n,d\geq 3$. In $1995$ Kollár \cite[Thm. 2]{KollarNonRatHyp1995} showed that a very general $X_d$ is not ruled if $n\geq 3$ and 
\begin{equation}\label{eq:bound_kollar}
d\geq 2\left \lceil\frac{n+3}{3}\right \rceil.
\end{equation}
Kollár's method is to degenerate $X_d$ to a $p:1$ cyclic cover $Y$ over a field of characteristic $p$. 
As $Y$ is necessarily singular he constructs a desingularization $\sigma: Y'\rightarrow Y$ and shows that, if the numerical condition (\ref{eq:bound_kollar}) holds, then $\wedge^{n-1} \Omega_{Y'}$ contains a big invertible sub-bundle $\cal L_{n,d}$. He then uses the positivity of $\cal L_{n,d}$ to show that $Y'$, and hence $Y$, is not ruled \cite[Lem. 7]{KollarNonRatHyp1995}.  As ruledness specializes by a result of Matsusaka \cite[App. 1.1]{Matsusaka}, he concludes that the complex hypersurface $X_d$ is not ruled.

Following Voisin's $2015$ breakthrough \cite{VoisinInventiones}, Totaro refined Kollár's method to show in \cite[Thm 2.1]{Totaro2016} that a very general hypersurface $X_d\subset \bb P_{\bb C}^{n+1}$ of degree $d$ is not stably rational if $n\geq 3$ and
\begin{equation}\label{eq:totaro}
d\geq 2\left \lceil\frac{n+2}{3}\right \rceil.
\end{equation}
Totaro makes use of the line bundle $\cal L_{n,d}$ of Kollár but only requires $\Gamma(\cal L_{n,d})\neq 0$, which explains the numerical improvement of (\ref{eq:totaro}) over (\ref{eq:bound_kollar}).
Successive works of Schreieder \cite{SchreiederAMS}, Moe--Ottem \cite{Moe2023}, Lange--Schreieder \cite{LangeSchreieder2024} and the very recent Hotchkiss \cite{Hotchkissquitic} and Böhning--von Bothmer--Pirutka--Zhang \cite{BohningBothmerPirutkaZhang2026} improve Totaro's bound and shows that for $d\geq 4$ a very general $X_d$ is not stably rational. \par
One can ask more generally 
what the minimal degree of a generically finite dominant rational map $X_d\dashrightarrow Z$ to a ruled variety $Z$ is. We denote this quantity by $irr_{rul}(X_d)$. If one additionally requires $Z$ to be rational one recovers the degree of irrationality $irr(X_d)$, 
first considered by Moh--Heinzer in \cite{MohHeinzer} and further studied by Bastianelli--De Poi--Ein--Lazarsfeld--Ullery in \cite{BDPELUcompositio}. Projecting from a point gives the inequality $d\geq irr(X_d)\geq irr_{rul}(X_d)$.

Chen--Stapleton prove in \cite[Lem. 2.3]{Chenstapleton} that if the line bundle $\cal L_{n,d}$ that appears in Kollár's construction separates $2r$ points in $Y'$ (see Definition \ref{def:separating}), then $irr_{rul}(Y)=irr_{rul}(Y')\geq \min\{p,r+1\}$, $p$ being the characteristic of the base field and $Y,Y'$ as above. 
Combining this observation with a degeneration result which shows that  $irr_{rul}$ can only drop upon specialization \cite[Thm. 1.1]{Chenstapleton}, they show that
\begin{equation}\label{eq:chenstap}
irr_{rul}(X_d)\geq \max_{p\leq d\,\,\text{prime}} \min\{p, \left \lceil \frac{(p+1)\lfloor \frac{d}{p} \rfloor -n}{2}\right \rceil\}
\end{equation}
for a very general hypersurface $X_d\subset \bb P^{n+1}_\bb C$ \cite[Prop. 2.5]{Chenstapleton}. Note that (\ref{eq:chenstap}) recovers (\ref{eq:bound_kollar}) by setting $irr_{rul}(X_d)\geq 2$ and $p=2$.
\par

Our main observation is that the line bundle $\cal L_{n,d}$ only needs to separate $r$ points in $Y'$ in order to prove Chen--Stapleton's inequality $irr_{rul}(X_d)\geq \min\{p,r+1\}$.
Since separating $r=1$ point is equivalent to having nontrivial global sections, this replaces the bigness condition on $\cal L_{n,d}$ for non-ruledness of $X_d$ with  $\Gamma(\cal L_{n,d})\neq 0$. 
Thus, as in Totaro’s approach to stable irrationality, non-ruledness follows from the existence of a nonzero global differential form. 
\begin{restatable}{theorem}{theoremain}\label{theorem:main}
Let $X_d\subset \bb P_{\bb C}^{n+1}$ be a very general hypersurface of degree $d$. Assume that $n\geq 3$ and $d\geq 2$, then 
\[
irr_{rul}(X_d)\geq\max_{p\leq d\,\,\text{prime}} \min\{p, (p+1)\left \lfloor \frac{d}{p} \right \rfloor -n\} .
\]
In particular $X_d$ is not ruled if
\begin{equation}\label{eq:nonrul}
d\geq 2\left \lceil\frac{n+2}{3}\right \rceil.
\end{equation}
\end{restatable}
The degree threshold (\ref{eq:nonrul}) improves  (\ref{eq:bound_kollar}) by two when $n=1\mod 3$ and agrees with it otherwise. 
We thus obtain that a very general quartic fourfold is not ruled. An alternative proof of this 
can also be given by degenerating a very general quartic fourfold to the union of a hyperplane and a very general cubic fourfold, whose non-ruledness was recently proven by Zhang \cite[Thm 1.3]{ZhangNonRuled} using Hodge atoms \cite{KKPY2025}.

As mentioned before very general cubic hypersurfaces are conjecturally not ruled
but one can see that these are birational to conic bundles. It is natural to ask if very general hypersurfaces $X_d$ with $d\geq 4$ and $n\geq 3$ might also admit birational conic bundle structures. Birational rigidity \cite{defernex}  settles this question negatively for $d=n+1$. Moreover, Kollár shows in  \cite[Thm 4.2]{KollarNonRatHyp1995} that if $n\geq 3$ and
\begin{equation}
d\geq 3\left \lceil\frac{n+3}{4}\right \rceil.
\end{equation}
 a very general $X_d$ is not birational to a conic bundle (see also \cite[Thm. V.5.14]{Kollar1996}). Analogously to $irr_{rul}$ one can consider the minimal degree of a rational generically finite dominant map $X_d\dashrightarrow Z$ to a variety $Z$ admitting a conic bundle structure. We denote this quantity by $irr_{con}(X_d)$.
Using the same method as for Theorem \ref{theorem:main} and a specialization result of Chen--Church--Ji--Stapleton \cite{chen_genus} we can bound $irr_{con}$. 
\begin{restatable}{theorem}{theoremain2}\label{theorem:main2}
Let $X_d\subset \bb P_{\bb C}^{n+1}$ be a very general hypersurface of degree $d$. Assume that $n\geq 3$ and $d\geq 3$, then 
\[
irr_{con}(X_d)\geq \max_{p\leq d\,\,\text{odd prime}} \min\{p, (p+1)\left \lfloor \frac{d}{p} \right \rfloor -n\}
\]
In particular $X_d$ is not birational to a conic bundle if
\begin{equation}\label{eq:nonconicbundle}
d\geq 3\left \lceil\frac{n+2}{4}\right \rceil.
\end{equation}
\end{restatable}

\paragraph{Acknowledgements}
I would like to thank Morten Lüders for introducing me to the subjects of the note and for many helpful discussions and suggestions.
Furthermore, I thank Nathan Chen, János Kollár and Jan Lange for
helpful comments. 
The author acknowledges support by Deutsche Forschungsgemeinschaft  (DFG) through the Collaborative Research Centre TRR 326 ``Geometry and Arithmetic of Uniformized Structures'', project number 444845124.

\paragraph{AI Disclosure}
During the preparation of this paper, ChatGPT Sol $5.6$ was used for preliminary discussions and literature searches. 
In particular, it helped the author find a reference to the Jacobson correspondence in \cite[Thm. 2.6]{Tziolas24}, which motivated Section \ref{sect:obstrcution} and Proposition \ref{proposition:criterium}. The writing of the manuscript below was carried out exclusively
by the author, who takes full responsibility for its content.

\section{An obstruction to maps to ruled varieties}\label{sect:obstrcution}
In this section we present a criterion to bound the degree of non-ruledness $irr_{rul}$. First recall what it means for sections of a vector bundle to separate points.
\begin{definition}\cite[Def. 2.1]{Chenstapleton}\label{def:separating}
Let $X$ be a variety over an algebraically closed field $k$ and let $\cal L$ be a vector bundle on $X$. We say that $\cal L$ separates $r$ points on an open $\emptyset \neq U\subset X$ if for any distinct $r$ points $x_1, ..., x_{r}\in U$ the evaluation $H^0(X,\cal L)\to \cal L_{x_1}\oplus ...\oplus \cal L_{x_r}$ is surjective.
\end{definition}
In particular if $\cal L$ is a line bundle it separates one point on some open if and only if there is  a nontrivial global section $0\neq s\in H^0(X,\cal L)$ and more generally it separates $r$ points $x_1, \ldots, x_{r}$ if and only if for each $1\leq i\leq r$ there is a section $a\in H^0(X,\cal L) $ that vanishes at $x_1,...,\hat{x_{i}},...x_r$ and doesn't vanish at $x_i$.
 Before coming to our criterion we first recall the result of Chen--Stapleton \cite{Chenstapleton}.
\begin{lemma}\label{lemma:chenstap}\cite[Lem. 2.3]{Chenstapleton}\cite[Lem. 7]{KollarNonRatHyp1995}
Let $X$ be a projective variety over an algebraically closed field and $\cal L$ a line bundle on $X$ such that $H^0(X,\cal L)$ separates $2r$ points on an open set. Suppose that there is an injection $\cal L\hookrightarrow \wedge^i\Omega_{X}$ for some $i>0$. Then every dominant, separable, generically finite rational map from $X$
to a ruled variety has degree at least $r+1$.
\end{lemma}
Our observation is that, in the special case where $\cal L$ is obtained from the pullback of differential forms on a smooth variety via an inseparable map, it suffices for $\cal L$ to separate just $r$ points to obtain the same conclusion as Lemma \ref{lemma:chenstap}.
For the rest of the section we work over an algebraically closed field $k$ of characteristic $p>0$.
\begin{proposition}\label{proposition:criterium}
    Let $\pi:X\to Z$ be a dominant purely inseparable morphism of degree $[k(X):k(Z)]=p$ with $X,Z$ smooth proper $n$-dimensional varieties.
    Fix $1\le i\le n-1$ and assume $H^0(Z, \Omega^{i}_Z)=0$.
    Set $\cal M:=im(\pi^*\Omega_Z\to \Omega_X)$. Let $\cal L$ be a vector bundle with an injection $\cal L\hookrightarrow \Omega_X^i$ such that $\cal L_\eta =\wedge^{i}\cal M_\eta$ for $\eta\in X$ the generic point. If $\cal L$ separates $r$ points on some open $\emptyset\neq U\subset X$, then every dominant, separable, generically finite rational map from $X$ to a ruled variety has degree at least $r+1$, in particular \[irr_{rul}(X)\geq \min\{p,r+1\}.\] 
    If $p\neq 2$ the same degree bound holds for targets $V$ admitting a dominant morphism $q:V\to Y$ which has generic fiber a smooth conic, that is $irr_{con}(X)\geq \min\{p,r+1\}$.
\end{proposition}
For what follows let $\phi:X\dashrightarrow Y\times \bb P^1$ be a rational dominant separable morphism of degree $b:=\deg(\phi)$. Denote the trace $Tr_{k(X)/k(Y)(t)}:\Omega^i_{k(X)/k}\to \Omega^i_{k(Y)(t)/k}$.
\begin{lemma}\label{lemma:lift}
    For any form $\alpha\in H^0(X,\Omega_X^i)$ there exists a unique form $\alpha'\in \Omega^i_{k(Y)/k}$ that restricts to its trace $Tr_{k(X)/k(Y)(t)}(\alpha)\in \Omega^i_{k(Y)(t)/k}$. Moreover for a rational dominant map $\psi:Z\dashrightarrow Y$ the pullback extends to a global form $\psi^*(\alpha')\in H^0(Z,\Omega^i_{Z/k})$.
\end{lemma}
\begin{proof}
Take the normalisation of the graph of $\phi$ with projections
\begin{equation}\label{diag:graohtrace}
\begin{tikzcd}
	& \Gamma & \\
	X && {Y\times \bb P^1}
	\arrow["{p_1}"', from=1-2, to=2-1]
	\arrow["{p_2}", from=1-2, to=2-3]
	\arrow["\phi", dashed, from=2-1, to=2-3]
\end{tikzcd}
\end{equation}
By the work of de Jong--Starr \cite[Prop. 3.3]{deJongStarr},\cite{Dejongerratum} there is a trace $Tr^i(p_1^*\alpha)\in H^0(Y\times \bb P^1, \Omega_{Y\times \bb P^1}^{i\vee\vee})$ which restricts to the trace $Tr_{k(X)/k(Y)(t)}(\alpha)$.  Moreover \[
\Omega_{Y\times \bb P^1}^{i\vee\vee}\cong pr_1^*\Omega_{Y}^{i\vee\vee}\oplus (pr_1^*\Omega_{Y}^{i-1\vee\vee}\otimes pr_2^*\Omega_{\bb P^1}^{1}).\]
As $ \Omega_{\bb P^1}$ does not have nontrivial global sections, any $i$-form must be pulled back from a form $\alpha'\in \Omega^i_{k(Y)/k}$. Uniqueness of $\alpha'$ follows from the injectivity of $\Omega^i_{k(Y)/k}\to \Omega^i_{k(Y)(t)/k}$. The same construction works with $Y$ replaced with any proper normal model $Y'$ of $k(Y)$, hence $\alpha'$ is regular at every codimension one point of $Y'$. Now let $\psi:Z\dashrightarrow Y$ be dominant, since $Z$ is smooth it is enough to show that $\psi^*\alpha'$ is regular at every codimension one point $z\in Z$.  Let $\nu_z$ be the associated valuation on $k(Z)$. If $\nu_z|_{k(Y)}$ is trivial then $k(Y)\subset \cal O_{Z,z}$ and $\psi^*\alpha'$ extends. If $\nu_z|_{k(Y)}$ is not trivial then it must be a divisorial valuation. Choose a normal model $Y'$ of $k(Y)$ on which the valuation $\nu_z|_{k(Y)}$ comes from a prime divisor $E\subset Y'$. By the previous paragraph the form $\alpha'$ is regular at the generic point of $E$. Hence $\psi^*\alpha'$ must be regular at $z$.
\end{proof}
\begin{remark}(Conic bundles)\label{remark:coniccase}
    Let $p\neq 2$ and $q:V\to Y$ be a dominant morphism of normal
varieties whose generic fibre is a smooth conic, and let $X\dashrightarrow V$ be dominant, generically finite and separable of degree $b$. Then if $char(k)\neq 2$ the conclusions of Lemma \ref{lemma:lift} still hold with $V$ instead of $Y\times \bb P^1$.
To prove this, apply the same construction as in (\ref{diag:graohtrace}) to $X\dashrightarrow V$. The trace is regular on the generic conic $V_{k(Y)}$ as it lies in the smooth locus of $V$. Taking global sections gives $H^0(V_{k(Y)}, \Omega^i_{V_{k(Y)}/k})\cong \Omega^i_{k(Y)/k}$. This proves the existence of the lift $\alpha'$ and uniqueness follows from the smoothness of $V_{k(Y)}$. Take a degree two morphism $V_{k(Y)}\to \bb P^1_{k(Y)}$. This gives a separable degree two extension $k(V)/k(Y)(t)$. Using Lemma \ref{lemma:lift} we have \[
Tr_{k(X)/k(Y)(t)}(\alpha)=Tr_{k(V)/k(Y)(t)}(\alpha')=2\alpha'\in \Omega_{k(Y)(t)/k}^i
\]
and since we assumed that $2$ is invertible the claim follows from Lemma \ref{lemma:lift}.
\end{remark}
As $\phi$ is separable it induces an isomorphism $k(X)\otimes_{k(Y)(t)} \Omega^i_{k(Y)(t)/k}\cong \Omega^i_{k(X)/k}$, hence the pullback $\phi^*:\Omega^i_{k(Y)(t)/k}\to \Omega^i_{k(X)/k}$ is injective and for $\alpha\in \Omega^i_{k(X)/k}$ we will identify $Tr_{k(X)/k(Y)(t)}(\alpha)$ and $\phi^*Tr_{k(X)/k(Y)(t)}(\alpha)$. 

\begin{lemma}\label{lemma:technical}
    Let $\cal L\hookrightarrow \Omega^i_{X/k}$ be a vector bundle. If $\cal L$ separates $b$ points on a nonempty open of $X$, then
    \begin{equation*}
        \cal L_\eta\subset \langle Tr_{k(X)/k(Y)(t)}(\alpha)\,|\,\alpha\in H^0(X,\cal L)\rangle_{k(X)}\,\,\,\text{in}\,\,\,\Omega^i_{k(X)/k}.
    \end{equation*}
\end{lemma}
\begin{proof}
As $\phi$ is separable there exist open affines $\Spec B\subset X$ and $\Spec A\subset Y\times \bb P^1$ where $\phi$ restricts to an étale map $A\rightarrow B$ of degree $b$.  For any closed point $y\in \Spec A$ we have $\phi^{-1}(y)=\{x_1,...,x_b\}\subset \Spec B(k)$. For any $0\neq \lambda\in \cal L_{x_1}$ separation gives a section $\alpha \in H^0(X, \cal L)$ with restrictions
\begin{equation*}\label{eq:retsr}
(\alpha|_{x_1},...,\alpha|_{x_b})=(\lambda,0,...,0)\in \cal L_{x_1}\oplus ...\oplus \cal L_{x_b}.
\end{equation*}
 The étale map $A\to B$ also induces a trace map for differential forms. We obtain the following commutative diagram
\begin{equation}\label{diag:normetale}
\begin{tikzcd}
	& {H^0(X,\cal L)} & \\
	{\Omega^i_{k(X)/k}} & {\Omega^i_{B/k}} & {\bigoplus_{j=1}^b\Omega^i_{B/k}\otimes_B\kappa(x_j)} \\
	{\Omega^{i}_{k(Y)(t)/k}} & {\Omega^{i}_{A/k}} & {\Omega^{i}_{A/k}\otimes_A\kappa(y)}
	\arrow[from=1-2, to=2-1]
	\arrow[from=1-2, to=2-2]
	\arrow["{\oplus\,\_|_{x_j}}", from=1-2, to=2-3]
	\arrow["{Tr_{k(X)/k(Y)(t)}}"', from=2-1, to=3-1]
	\arrow[hook', from=2-2, to=2-1]
	\arrow["{\oplus\,\_|_{x_j}}", from=2-2, to=2-3]
	\arrow["{Tr_{B/A}}"', from=2-2, to=3-2]
	\arrow["\rho"', from=2-3, to=3-3]
	\arrow[hook', from=3-2, to=3-1]
	\arrow["{|_y}", from=3-2, to=3-3]
\end{tikzcd}
\end{equation}
where the morphism $\rho$ is defined using the isomorphism $\Omega^{i}_{B/k}\otimes_B\kappa(x_j)\cong \Omega^{i}_{A/k}\otimes _A\kappa (y)$ (which follows from the fact that $k$ is algebraically closed and $A\to B$ is étale) and the natural sum morphism $\bigoplus \Omega^{i}_{A/k}\otimes _A\kappa (y)\to \Omega^{i}_{A/k}\otimes _A\kappa (y)$. Consequently \[
\phi^*Tr_{B/A}(\alpha)|_{x_1}=\lambda\in \Omega^i_{B/k}\otimes_B \kappa (x_1).
\]
Consider the $B$-linear map 
\begin{equation}\label{eq:technicalmap}
f: H^0(\Spec B, \cal L)\to \frac{\Omega^i_{B/k}}{\langle Tr_{B/A}(\alpha)\,|\,\alpha\in H^0(X,\cal L)\rangle_{B}}
\end{equation}
we have shown that $f$ vanishes restricted to any points $x_1\in \Spec B$, hence it also vanishes upon shrinking to an open where the target of $f$ is locally free. As the trace commutes with localisation the restriction of (\ref{eq:technicalmap}) to the generic point $\eta$ of $X$ proves the claim.    
\end{proof}
The main step in the proof of Proposition \ref{proposition:criterium} is to show that a rational dominant separable morphism $\phi:X\dashrightarrow Y\times \bb P^1$ of degree $b:=\deg(\phi)\leq r$ factorizes over $Z$, meaning that there is a rational map $\psi$ that makes the square 
\begin{equation}\label{eq:square}
\begin{tikzcd}
	X & {Y\times \bb P^1} \\
	Z & Y
	\arrow["\phi", dashed, from=1-1, to=1-2]
	\arrow["\pi"', from=1-1, to=2-1]
	\arrow["{pr_1}", from=1-2, to=2-2]
	\arrow["{\exists \psi}", dashed, from=2-1, to=2-2]
\end{tikzcd}
\end{equation}
commute. The factorization will be given by Lemmas \ref{lemma:lift} and \ref{lemma:technical} and a simple consequence of the Jacobson correspondence \cite[Thm 2.6]{Tziolas24}, which we recall next. For a $k$-algebra $A$ we denote $Der_k(A)$ the $k$-vector space of $k$-linear derivations $d:A\to A$.

\begin{proposition}(Jacobson correspondence)\cite[Thm 2.6]{Tziolas24}\label{prop:jacobson}
    Let $L/k$ be a finitely generated field extension. There is an inclusion-reversing bijection between fields $K$ with a $k$-algebra structure and inclusions $L^p\subset K\subset L$ and the set of restricted Lie subalgebras $S\subset Der_k(L)$ which are finite-dimensional $L$-vector spaces. The mappings are given by $K\mapsto Der_K(L)$ and $S\mapsto L^{S}$, where $L^S$ denotes the intersection of the kernels of the derivations in $S$. 
\end{proposition}

\begin{proof}[Proof of Proposition \ref{proposition:criterium}]
Let $\phi:X\dashrightarrow Y\times \bb P^1$ be a rational dominant separable morphism of degree $b:=\deg(\phi)\leq r$. We show that if $\cal L$ separates $b$ points in a nonempty open of $X$ then $Z$ must have nonzero global $i$-forms, which contradicts the assumptions. Consider the two $k(X)$ -vector subspaces of $\Omega^1_{k(X)/k}$ given by
\begin{equation}
    W_Z:=im(k(X)\otimes \Omega^1_{k(Z)/k}\to \Omega^1_{k(X)/k}),\,\,\,W_Y:=im(k(X)\otimes \Omega^1_{k(Y)/k}\to \Omega^1_{k(X)/k}).
\end{equation}
By assumption $\bigwedge^iW_Z=\cal L_\eta$, moreover we have 
\begin{equation}
    \cal L_\eta\stackrel{Lem. \ref{lemma:technical}}{\subset} \langle Tr_{k(X)/k(Y)(t)}(\alpha)\,|\,\alpha\in H^0(X,\cal L)\rangle_{k(X)}\stackrel{Lem. \ref{lemma:lift}}{\subset} \bigwedge^i W_Y
\end{equation}
Both $W_Y$ and $W_Z$ have the same dimension $n-1$ so $\bigwedge^iW_Z\subset \bigwedge^iW_Y$ directly implies $W_Z=W_Y$. It follows that $Der_{k(Y)}(k(X))=Der_{k(Z)}(k(X))$ and the Jacobson correspondence (Proposition \ref{prop:jacobson}) therefore gives 
\[
k(Y)\subset k(X)^{Der_{k(Y)}(k(X))}=k(X)^{Der_{k(Z)}(k(X))}=k(Z).
\]
We have the following inclusions of $k$-algebras
\begin{equation}\label{diamond}
    \begin{tikzcd}
	& {k(X)} & \\
	{k(Z)} && {k(Y)(t)} \\
	& {k(Y)}
	\arrow["\pi", hook, from=2-1, to=1-2]
	\arrow["\phi"', hook', from=2-3, to=1-2]
	\arrow["\psi", hook, from=3-2, to=2-1]
	\arrow["\iota"', hook', from=3-2, to=2-3]
\end{tikzcd}
\end{equation}
As $\cal L_\eta$ is not zero by Lemma \ref{lemma:technical} we can find a form $\alpha\in H^0(X,\Omega^i_{X/k})$ such that $Tr_{k(X)/k(Y)(t)}(\alpha)\neq 0$. By Lemma \ref{lemma:lift} there exists a form $\alpha'\in \Omega^i_{k(Y)/k}$ that restricts to the trace of $\alpha$, using the commutativity of (\ref{diamond}) we see 
\[
\pi^*\psi^*\alpha'=\phi^*\iota^*\alpha'=Tr_{k(X)/k(Y)(t)}(\alpha)\neq0
\]
But by Lemma \ref{lemma:lift} the form $\psi^*\alpha'$ extends to a global differential in $\Omega_{Z/k}^i$ and the proof of the first claim is concluded. For the conic-bundle assertion, replace $k(Y)(t)$ throughout with $k(V)$. The statement of Lemma \ref{lemma:technical} translates verbatim and by Remark \ref{remark:coniccase} Lemma \ref{lemma:lift} applies as well. The same arguments as in the previous paragraph show the stated degree bound.
\end{proof}

\section{Proofs of Theorems \ref{theorem:main} and \ref{theorem:main2}}
The proofs follow closely those of \cite[Thm. A]{Chenstapleton} of Chen--Stapleton and \cite[Thm. A]{chen_genus} of Chen--Church--Ji--Stapleton, we include sketches for completeness.

We first recall the common degeneration construction.
Start with a very general hypersurface $X_d\subset \bb P^{n+1}_\bb C$. For a prime $d\geq p\geq 2$ write $d=pe+\lambda$ for $0\leq  \lambda< p$. Let $k$ be an algebraically closed field of characteristic $p$. Using a construction given by Mori in his thesis \cite[4.3]{Mori} we can degenerate $X_d$ to a union of $\lambda$ hyperplanes and a very general cyclic $p:1$ cover $\pi:X\to Z$ of a hypersurface $Z\subset \bb P_k^{n+1}$ of degree $e$ branched over a degree $pe$ hypersurface in $Z$ (see also \cite[Constr. 5]{KollarNonRatHyp1995}). The map $\pi:X\to Z$ is inseparable and of degree $p$. By the assumption $n\geq 3$ Kollár shows in \cite{KollarNonRatHyp1995} that there is an explicit resolution of singularities $\rho:X'\to X$ with a line bundle
\[
\cal L_{n,d}:=(\pi\rho)^*(\omega_Z\otimes \cal O_Z(pe))\cong ( \pi\rho)^*\cal O_Z((p+1)e-n-2)
\]
which embeds in the $n-1$-differentials of $X'$, $\cal L_{n,d}\hookrightarrow \Omega^{n-1}_{X'}$ and coincides with $\wedge^{n-1}im((\pi\rho)^*\Omega_Z\to \Omega_{X'})$ at the generic point of $X'$. When $(p+1)e-n-2\geq 0$ the line bundle $\cal L_{n,d}$ separates $(p+1)e-n-1$ points in $X'$. Moreover by weak Lefschetz \cite[Def. 5.1, Prop. 5.3]{AntieuBhattmathew}
\[
H^0(\bb P_k^{n+1}, \Omega^{n-1}_{\bb P_k^{n+1}})\xrightarrow{\cong}H^0(Z, \Omega_Z^{n-1})
\]
hence $H^0(Z, \Omega_Z^{n-1})=0$. 
\begin{proof}[Proof of Theorem \ref{theorem:main}]
Fix a prime $p$ with $p\leq d$. We can assume that $(p+1)e-n-2\geq 0$.
We can apply Proposition \ref{proposition:criterium} to the composition $\pi\rho:X'\to Z$ with $i=n-1$ and $r=(p+1)e-n-1$ and obtain that 
\[
irr_{rul}(X)=irr_{rul}(X')\geq \min\{p,(p+1)e-n\}.
\]
 Using the fact that $irr_{rul}$ can only decrease in the degeneration, proved by Chen--Stapleton in \cite[Thm. 1.1]{Chenstapleton}, we conclude that \begin{equation}\label{eq:bound}
irr_{rul}(X_d)\geq \min\{p,(p+1)e-n\}=\min\{p,(p+1)\left \lfloor \frac{d}{p}\right \rfloor-n\}
    \end{equation}
Finally setting $p=2$ we see that (\ref{eq:bound}) provides the non-ruledness of $X_d$ if $d\geq 2\lceil\frac{n+2}{3}\rceil$.
\end{proof}
\begin{lemma}\label{lemma:secondthm}
    If $p\neq 2$ then in the degeneration above we have $irr_{con}(X_d)\geq irr_{con}(X)$.
\end{lemma}
\begin{proof}
We apply the construction of \cite[Thm 1.1]{Chenstapleton}. 
Start with a degree $b$ map from $X_d$ to a conic bundle $X_d\dashrightarrow Y$.
After a finite base change defining both the map and the conic
fibration, choose a normal integral scheme $\cal Y$ which is proper and flat over the DVR and has geometric generic fiber $Y$. Likewise choose normal integral model $\cal X$ for the degenerating family of $X_d$ to $X$. Note  that $X$ has multiplicity one hence it is birational to a component of the special fiber of $\cal X$.  By an argument of Abhyankar--Zariski \cite[Lem. 2.22]{Kollarsingularieties} we may blow up and normalize $\cal Y$ so that the generic point of $X$ maps
to the generic point of a special-fiber component of $\cal Y$. By \cite[Tag~09E5]{stacks-project} the induced dominant rational map from $X$ to this component
has degree at most $b$. By \cite[Prop.~3.6(1)]{chen_genus} we have that the component of $\cal Y$ to which $X$ maps has fibering genus equal to zero.
After normalization and Stein factorization, its generic curve is regular, proper, and of genus zero.  Since $p\neq 2$ it is a smooth conic by \cite[Thm. 2.8]{chen_genus}. Taking
$b=irr_{con}(X_d)$ proves the assertion.
\end{proof}
\begin{proof}[Proof of Theorem \ref{theorem:main2}]
The application of Proposition \ref{proposition:criterium} in the proof of Theorem \ref{theorem:main} repeats verbatim and we conclude using Lemma \ref{lemma:secondthm}. Setting $p=3$ provides the bound for $X_d$ not being birational to a conic bundle.
\end{proof}

\bibliographystyle{acm}
\bibliography{bib}

@misc{BohningBothmerPirutkaZhang2026,
  author        = {Christian B{\"o}hning and Hans-Christian Graf von Bothmer and Alena Pirutka and Zhijia Zhang},
  title         = {Stable irrationality of quartic hypersurfaces},
  year          = {2026},
  eprint        = {2611.15702},
  archivePrefix = {arXiv},
  primaryClass  = {math.AG},
  doi           = {10.48550/arXiv.2611.15702},
  url           = {https://arxiv.org/abs/2611.15702}
}

@misc{ZhangNonRuled,
  author = {Zhang, Nantao},
  title = {Non-ruledness via Hodge atoms},
  note = {Unpublished draft},
  url = {https://zenith-john.github.io/file/Non-ruled.pdf}
}

@article{AntieuBhattmathew,
 author = {Antieau, Benjamin and Bhatt, Bhargav and Mathew, Akhil},
 title = {Counterexamples to {Hochschild}-{Kostant}-{Rosenberg} in characteristic {{\(p\)}}},
 fjournal = {Forum of Mathematics, Sigma},
 journal = {Forum Math. Sigma},
 issn = {2050-5094},
 volume = {9},
 pages = {26},
 note = {Id/No e49},
 year = {2021},
 language = {English},
 doi = {10.1017/fms.2021.20},
 zbMATH = {7377623},
 Zbl = {1477.13032}
}

@article{BDPELUcompositio,
 author = {Bastianelli, Francesco and De Poi, Pietro and Ein, Lawrence and Lazarsfeld, Robert and Ullery, Brooke},
 title = {Measures of irrationality for hypersurfaces of large degree},
 fjournal = {Compositio Mathematica},
 journal = {Compos. Math.},
 issn = {0010-437X},
 volume = {153},
 number = {11},
 pages = {2368--2393},
 year = {2017},
 language = {English},
 doi = {10.1112/S0010437X17007436},
 zbMATH = {6810456},
 Zbl = {1387.14118}
}

@misc{Hotchkissquitic,
 author = {James Hotchkiss},
 title = {The rationality problem for hypersurfaces of degree at least five},
 year = {2026},
 howpublished = {Preprint, {arXiv}:2610.10509 [math.{AG}] (2026)},
 url = {https://arxiv.org/abs/2610.10509},
 arXiv = {arXiv:2610.10509}
}

@book{Kollarsingularieties,
 author = {Koll{\'a}r, J{\'a}nos},
 title = {Singularities of the minimal model program. {With} the collaboration of {S{\'a}ndor} {Kov{\'a}cs}},
 fseries = {Cambridge Tracts in Mathematics},
 series = {Camb. Tracts Math.},
 issn = {0950-6284},
 volume = {200},
 isbn = {978-1-107-03534-8; 978-1-139-54789-5},
 year = {2013},
 publisher = {Cambridge: Cambridge University Press},
 language = {English},
 doi = {10.1017/CBO9781139547895.002},
 zbMATH = {6148846},
 Zbl = {1282.14028}
}

@article{chen_genus,
  author        = {Chen, Nathan and Church, Benjamin and
                   Ji, Lena and Stapleton, David},
  title         = {The fibering genus of {Fano} hypersurfaces},
  journal       = {Algebraic Geometry and Physics},
  note          = {To appear},
  eprint        = {2308.12401},
  archivePrefix = {arXiv},
  primaryClass  = {math.AG},
  url           = {https://arxiv.org/abs/2308.12401}
}

@article{Dejongerratum,
 author = {de Jong, A. J. and Starr, Jason},
 title = {Erratum: {Cubic} fourfolds and spaces of rational curves},
 fjournal = {Illinois Journal of Mathematics},
 journal = {Ill. J. Math.},
 issn = {0019-2082},
 volume = {52},
 number = {1},
 pages = {345--346},
 year = {2008},
 language = {English},
 zbMATH = {5567376},
 Zbl = {1162.14302}
}

@article{MohHeinzer,
 author = {Moh, T. T. and Heinzer, W.},
 title = {On the {L{\"u}roth} semigroup and {Weierstrass} canonical divisors},
 fjournal = {Journal of Algebra},
 journal = {J. Algebra},
 issn = {0021-8693},
 volume = {77},
 pages = {62--73},
 year = {1982},
 language = {English},
 doi = {10.1016/0021-8693(82)90277-0},
 zbMATH = {3845738},
 Zbl = {0533.14010}
}

@article{deJongStarr,
  author  = {de Jong, A. J. and Starr, Jason},
  title   = {Cubic fourfolds and spaces of rational curves},
  journal = {Illinois Journal of Mathematics},
  volume  = {48},
  number  = {2},
  year    = {2004},
  pages   = {415--450},
  doi     = {10.1215/ijm/1258138390}
}

@article{Mori,
 author = {Mori, Shigefumi},
 title = {On a generalization of complete intersections},
 fjournal = {Journal of Mathematics of Kyoto University},
 journal = {J. Math. Kyoto Univ.},
 issn = {0023-608X},
 volume = {15},
 pages = {619--646},
 year = {1975},
 language = {English},
 doi = {10.1215/kjm/1250523007},
 zbMATH = {3518949},
 Zbl = {0332.14019}
}

@article{Matsusaka,
 author = {Matsusaka, T.},
 title = {Algebraic deformations of polarized varieties},
 fjournal = {Nagoya Mathematical Journal},
 journal = {Nagoya Math. J.},
 issn = {0027-7630},
 volume = {31},
 pages = {185--245},
 year = {1968},
 language = {English},
 doi = {10.1017/S0027763000012733},
 zbMATH = {3268748},
 Zbl = {0167.49501}
}

@article{defernex,
 author = {de Fernex, Tommaso},
 title = {Birationally rigid hypersurfaces},
 fjournal = {Inventiones Mathematicae},
 journal = {Invent. Math.},
 issn = {0020-9910},
 volume = {192},
 number = {3},
 pages = {533--566},
 year = {2013},
 language = {English},
 doi = {10.1007/s00222-012-0417-0},
 zbMATH = {6176666},
 Zbl = {1279.14019}
}

@article{Chenstapleton,
 author = {Chen, Nathan and Stapleton, David},
 title = {Fano hypersurfaces with arbitrarily large degrees of irrationality},
 fjournal = {Forum of Mathematics, Sigma},
 journal = {Forum Math. Sigma},
 issn = {2050-5094},
 volume = {8},
 pages = {12},
 note = {Id/No e24},
 year = {2020},
 language = {English},
 doi = {10.1017/fms.2020.20},
 zbMATH = {7205611},
 Zbl = {1454.14034}
}

@article{Tziolas24,
 author = {Tziolas, Nikolaos},
 title = {Topics in group schemes and surfaces in positive characteristic},
 fjournal = {Annali dell'Universit{\`a} di Ferrara. Sezione VII. Scienze Matematiche},
 journal = {Ann. Univ. Ferrara, Sez. VII, Sci. Mat.},
 issn = {0430-3202},
 volume = {70},
 number = {3},
 pages = {891--954},
 year = {2024},
 language = {English},
 doi = {10.1007/s11565-024-00514-8},
 zbMATH = {7901292},
 Zbl = {1569.14051}
}

@article{VoisinInventiones,
 author = {Voisin, Claire},
 title = {Unirational threefolds with no universal codimension {{\(2\)}} cycle},
 fjournal = {Inventiones Mathematicae},
 journal = {Invent. Math.},
 issn = {0020-9910},
 volume = {201},
 number = {1},
 pages = {207--237},
 year = {2015},
 language = {English},
 doi = {10.1007/s00222-014-0551-y},
 zbMATH = {6468725},
 Zbl = {1327.14223}
}

@misc{KKPY2025,
      title={Birational Invariants from Hodge Structures and Quantum Multiplication}, 
      author={Ludmil Katzarkov and Maxim Kontsevich and Tony Pantev and Tony Yue YU},
      year={2025},
      eprint={2508.05105},
      archivePrefix={arXiv},
      primaryClass={math.AG},
      url={https://arxiv.org/abs/2508.05105}, 
}

@article{KollarNonRatHyp1995,
 author = {Koll{\'a}r, J{\'a}nos},
 title = {Nonrational hypersurfaces},
 fjournal = {Journal of the American Mathematical Society},
 journal = {J. Am. Math. Soc.},
 issn = {0894-0347},
 volume = {8},
 number = {1},
 pages = {241--249},
 year = {1995},
 language = {English},
 doi = {10.2307/2152888},
 zbMATH = {788746},
 Zbl = {0839.14031}
}

@article{Totaro2016,
 author = {Totaro, Burt},
 title = {Hypersurfaces that are not stably rational},
 fjournal = {Journal of the American Mathematical Society},
 journal = {J. Am. Math. Soc.},
 issn = {0894-0347},
 volume = {29},
 number = {3},
 pages = {883--891},
 year = {2016},
 language = {English},
 doi = {10.1090/jams/840},
 zbMATH = {6572969},
 Zbl = {1376.14017}
}

@article{SchreiederAMS,
 author = {Schreieder, Stefan},
 title = {Stably irrational hypersurfaces of small slopes},
 fjournal = {Journal of the American Mathematical Society},
 journal = {J. Am. Math. Soc.},
 issn = {0894-0347},
 volume = {32},
 number = {4},
 pages = {1171--1199},
 year = {2019},
 language = {English},
 doi = {10.1090/jams/928},
 zbMATH = {7121119},
 Zbl = {1442.14138}
}

@article{LangeSchreieder2024,
 author = {Lange, Jan and Schreieder, Stefan},
 title = {On the rationality problem for low degree hypersurfaces},
 fjournal = {Forum of Mathematics, Pi},
 journal = {Forum Math. Pi},
 issn = {2050-5086},
 volume = {13},
 pages = {31},
 note = {Id/No e25},
 year = {2025},
 language = {English},
 doi = {10.1017/fmp.2025.10016},
 zbMATH = {8141826}
}

@misc{Moe2023,
 author = {Moe, Simen Westbye},
 title = {On {Stable} {Rationality} of {Polytopes}},
 year = {2023},
 howpublished = {Preprint, {arXiv}:2311.01144 [math.{AG}] (2023)},
 url = {https://arxiv.org/abs/2311.01144},
 arXiv = {arXiv:2311.01144}
}

@book {Kollar1996,
    AUTHOR = {Koll\'{a}r, J\'{a}nos},
     TITLE = {Rational curves on algebraic varieties},
    SERIES = {Ergebnisse der Mathematik und ihrer Grenzgebiete. 3. Folge. A
              Series of Modern Surveys in Mathematics [Results in
              Mathematics and Related Areas. 3rd Series. A Series of Modern
              Surveys in Mathematics]},
    VOLUME = {32},
 PUBLISHER = {Springer-Verlag, Berlin},
      YEAR = {1996},
     PAGES = {viii+320},
      ISBN = {3-540-60168-6},
   MRCLASS = {14-02 (14C05 14E05 14F17 14J45)},
  MRNUMBER = {1440180},
MRREVIEWER = {Yuri G. Prokhorov},
       DOI = {10.1007/978-3-662-03276-3},
       URL = {https://doi.org/10.1007/978-3-662-03276-3},
}

@misc{stacks-project,
    shorthand    = {Stacks},
    author       = {{The Stacks Project Authors}},
    title        = {\textit{Stacks Project}},
    howpublished = {\url{https://stacks.math.columbia.edu}},
    year         = {2018},
  }

\end{document}